\documentstyle[amsfonts,amscd,amssymb,latexsym,12pt]{article}

\title {\huge \bf A Principle for Critical Point under Generalized Regular Constraint  and Ill- Posed Lagrange
Multipliers under Non-Regular Constraints  }
\author{Ma Jipu$^{1,2}$  }
\date{}

\begin{document}
\maketitle

\newtheorem{th}{Theorem}[section]
\newtheorem{ex}{Example}
\newtheorem{re}{Remark}[section]
\newtheorem{lem}{Lemma}[section]
\newtheorem{de}{Definition}[section]
\newtheorem{co}{Corollary}[section]
\renewcommand{\theequation}{\thesection.\arabic{equation}}
\begin{quote}
 ABSTRACT. \ In this paper, a kind of non regular constraints and
  a principle for seeking critical point under the constraint are presented
  , where  no
  Lagrange multiplier is involved.  Let $E, F$ be two Banach spaces, $g: E\rightarrow F$ a $c^1$
   map defined on an open set $U$ in $E,$  and
 the constraint $S=$ the preimage $g^{-1}
(y_0), y_0\in F.$
 A main deference between  the non regular constraint and  regular constraint
  is that $g'(x)$ at any $x\in S$ is not surjective. Recently, the critical
  point theory under the non regular constraint is a concerned
   focus in optimization theory. The principle also suits the case
   of  regular constraint. Coordinately, the generalized regular
   constraint is introduced, and the critical point principle on generalized regular
   constraint is established. Let  $f: U \rightarrow \mathbb{R}$ be a
    nonlinear functional.
While the Lagrange multiplier $L$ in classical critical point
principle is considered, and its expression is given by using
generalized inverse ${g'}^+(x)$ of $g'(x)$ as follows : if $x\in S$
is a critical point of $f|_S,$ then $L=f'(x)\circ {g'}^+(x) \in F^*.
$ Moreover, it is proved that if $S$ is a regular constraint, then
the Lagrange multiplier $L$ is unique; otherwise, $L$ is ill-posed.
Hence, in case of the non regular constraint, it is very difficult
to solve Euler equations, however, it is often the case in
optimization theory.  So the principle here seems to be new and
applicable.
    By the way, the following theorem is proved; if $A\in B(E,F)$ is
double split, then the set of all generalized inverses of $A,$
$GI(A)$ is smooth
      diffeomorphic to certain Banach space. This is a new and
      interesting result in generalized inverse analysis.


\end{quote}
\vspace{1cm}
\section{ INTRODUCTION AND PRELIMINARY}\label{s1}
 \vspace{0.3cm}
 Let $E, F$ be two Banach spaces and $g:E\rightarrow F$ be a $c^1$  map defined on an
  open set $U$ in Banach space  $E.$
  Recall that a point $x\in U$ is said to be a regular point (or submersion ) of
 $g$ if
the Frech$\acute{e}$t differential $g'(x)$ is surjective and
$N(g'(x))$ splits  $E$, and $y_0\in F$ is a regular value of $g$
 provided
 the preimage $g^{-1}(y_0)$ is empty or consists only of regular points. It is known that the fundamental
 theorem on critical point theory under regular constraint holds:\\

\textbf{Theorem} (Preimage) \ If $y_0\in F$ is a regular value of
$c^1$ map $g:U\rightarrow F,$ then the preimage $S=g^{-1}(y_0)$ is a
$c^1$ submanifold of $U$ with the tangent space $T_xS=N(g'(x))$ for
any $x\in S.$
 Where $N(,)$ denotes the null space of the operator in the parenthesis. (For the details see \cite{abr}
  and  \cite{zei} .)                                                                                                 \\

 Recently, the concept of regular point has been extended to generalized regular point, i.e.,
  $x_0\in U$ is said to be a
generalized regular point of $g$ provided $g'(x_0)$ is double split,
 $F=R(g'(x_0))\oplus N_+,$ and there exists a neighborhood
$U_0\subset U$ of $x_0$ such that $R(g'(x)) \cap N_+=\{0\}$ for any
$x\in U_0,$ where $R(.)$ denotes the range of the operator in
 the parentheses. Obviously, regular point and immersion both are  generalized regular
points; when the rank of $g'(x_0), Rank(g'(x_0))< \infty,$ $x_0$ is
generalized regular point if and only if $x_0$ is subimmersion. It
is natural to define the generalized value $y_0\in F$ of $g$ as
$S=g^{-1}(y_0)$ is empty or consists only of generalized regular
points. We also have the following fundamental theorem of critical
point theory under generalized
regular constraint:\\

\textbf{Theorem} (Generalized Preimage) \ If $y_0\in F$ is
generalized regular value of $c^1$ map $g,$ then $S=g^{-1}(y_0)$ is
a $c^1$ submanifold of $U$ with the tangent space $T_xS=N(g'(x))$
for any $x\in S.$ (For the details see \cite{m6} and \cite {abr}.)\\

 Fortunately, we have the following complete rank theorem in advanced calculus
  to answer the locally conjugate problem proposed by Berger, M. in \cite {b},
 \\
 \textbf{Theorem} (Rank) \ Suppose that $g:U\subset E\rightarrow F$
is a $c^1$ map, $g'(x_0))$ is double split, and $g(x_0)=y_0, x_0\in
U.$ The following conclusion holds: there exist two neighborhoods
$U_0$ at $x_0,$ $V_0$ at 0, two local diffeomorphisms
$\varphi:U_0\rightarrow
\varphi (U_0)$ and $\psi: V_0\rightarrow \psi(V_0),$ such that\\
$$\varphi (x_0)=0, \ \varphi '(x_0)=I, \ \psi (0)=y_0, \ \psi
'(0)=I,$$and $$g(x)=(\psi \circ g'(x_0)\circ \phi )(x)$$for all
$x\in U_0$ if and only if $x_0$ is a generalized regular point of
$g.$ (For details see \cite {m3}, \cite {m7}, \cite {m8}, \cite {b},  \cite{abr}and \cite{zei}.)\\
(Note that the question on  rank theorem in advanced calculus
initially is
to find a sufficient condition such that the conclusion of the rank theorem above holds.)\\
 \indent There are many
equivalent conditions for generalized regular points, which are
convenient for analysis calculus.( For the details see
\cite{m1}, \cite{m3}, \cite{m9} and \cite{hm}.)\\
\indent By Theorem (Rank),the first main result in this paper, a
principle for seeking
 critical point under generalized regular constraint is given, in which no
 Lagrange multiplier is involved. Also, we present some  simple
 examples to illustrate its application and significance. Let
 $B(E,F)$ be the set of all linear bounded operators from $E$
 into
 $F;$ $A\in B(E,F)$ is called double split provided $R(A)$is
 closed and there exist closed subspaces $R_+\subset E$ and $N_+\subset F$
 such that $E=N(A)\oplus R_+$ and $F=R(A)\oplus N_+,$ where $N(A)$
 is the null space of $A.$ It is known well that there exists a
 generalized inverse $A^+\in B(F,E)$ for any double split operator
 $A,$ such that $AA^+A=A$ and $A^+AA^+=A^+.$ When $E$ and $F$ both are
 Hilbert spaces, a generalized inverse $A^+$ is said to be $M-P$
 inverse of $A$ provided $AA^+$ and $A^+A$ both are self adjoint.
 It is also known that the $M-P$ inverse of $A$ is unique, (see \cite{n}). Let
 ${q'}^+(x)$ be a generalized inverse of $g'(x).$  Then the second
 result is given as follows,suppose that $S$ is generalized regular
 constraint, and $f:U\rightarrow \mathbb{R}$ is a non linear functional;
  if $x$ is critical point of $f|_S,$
 then\\

 \quad \quad  $f'(x)-f'(x)\circ {q'}^+(x)\circ  q'(x)=0 $\quad  and  \quad  $ g(x)=y_0
 .$\\

 This shows that the Lagrange multiplier $L$ is a bounded linear functional
 on $F$ and $=f'(x)\circ{q'}^+(x)$ at the critical point $x$ of $f|_S;$
 moreover, when $S$ is a regular constraint, so that $R(g'(x))=F$ for
 any $x\in S,$  the Lagrange multiplier $L$ is unique since $f'(x)e=f'(x)$
$\circ {g'}^+(x)\circ g'(x)e \ \forall e\in E$ at the critical point
$x$ of $f|_S.$
  Specially, when $E$ and $F$ both are  Hilbert spaces, ${g'}^+(x)$
  can be a $M.-P.$ inverse of $g'(x),$ which is unique for any $x\in S.$
 Finally, let $A\in B(E,F)$ be double split, and $GI(A)$ the set of
all generalized inverses of $A$. The following theorem for
generalized inverse analysis is proved: $GI(A)$ is smooth
diffeomorphic to some Banach
 space, and then, using some results and presented method we can
 prove that the Lagrange multiplier $L$  under  non regular constraint
 is ill-posed. From this one can observe that it is very difficult
  to solve Euler equations in the case of non regular constraint.
 Hence besides the foundation for critical point under generalized regular
 constraint  the principle presented here for seeking critical points is a new and
 available way  since the Lagrange multiplier is no longer involved in
the principle.

\section{A Principle for Seeking Critical Point  }\label{s2}
 \vspace{0.3cm}
Let $g:U\rightarrow F$ be a $c^1$ map, and $y_0\in F$ a generalized
regular value of $g.$ By Theorem (Generalized Preimage), $S=$ the
preimage $g^{-1}(y_0)$ is a $c^1$ submanifold of $F.$ In what
follows, $S$ will be said to be the generalized regular constraint.
 Let $f$ be a non linear functional on $U.$ In this section, we discuss
 the critical point of $f$ under generalized regular constraint, and give a principle
  for seeking the critical point, while no Lagrange multiple is involved.

\begin{th}\label{2.1} If $x\in U$ is a critical point of $f|_S,$ where $S=g^{-1}(y_0),$
and $y_0\in F$ is a generalized regular value of $g,$ then

 \begin{equation}\label{2.1}
N(f'(x))\supset N(g'(x)) \quad {\rm{and}} \quad g(x)=y_0.
\end{equation}

\end{th}
{\it Proof.}\quad  Suppose that $x_0\in S$ is critical point of
$f|_S.$ Since $y_0$ is a generalized regular value, $x_0$ is a
generalized regular point of $g.$ By Theorem (Rank) there exist
neighborhoods $U_0$ at $x_0,$ $V_0$ at 0, local diffeomorphisms
$\varphi:U_0\rightarrow
\varphi (U_0)$ and $\psi: V_0\rightarrow \psi(V_0),$ such that\\
$$\varphi(x_0)=0, \ \varphi '(x_0)=I, \ \psi (0)=y_0, \ \psi '(0)=I,$$and
$$g(x)=(\psi \circ g'(x_0)\circ \varphi )(x) \quad \quad \forall x\in
U_0.$$ Since $\varphi :U_0\rightarrow \varphi (U_0)$ is a
diffeomorphism,  and $\varphi (x_0)=0,$ it is clear that there
exists a positive number $\varepsilon _0$ such that the following
relation for arbitrary fixed $h\in N(g'(x_0))$ holds : $x(t)=\varphi
^{-1}(th)\subset U_0$ for $|t|<\varepsilon_0.$ Then it follows
$$g(x(t))=(\psi \circ g'(x_0)\circ \varphi )(\varphi
^{-1}(th))=(\psi \circ g'(x_0))(th)=\psi (0)=y_0$$ (note  $\psi
(0)=y_0$), this shows the curve $x(t), -\varepsilon _0<t<\varepsilon
_0$ lies in $S,$ and so, $t=0$ is critical point of the function
$f(x(t)).$ Therefore,$$0=\frac{df(x(t))}{dt} |_{t=0}=
f'(x_0){(\varphi
 ^{-1}})'(0)h=f'(x_0)h \quad  \forall h\in N(g'(x_0)).$$ This proves $N(f'(x_0))\supset
 N(g'(x_0)).$ \hfill $\Box$\\
Specially, when $E$ is a Hilbert space, we have

\begin{th}\label{2.2} Suppose that $E$ is a Hilbert space, $S$ is a generalized regular
constraint and $f$ is a $c^1$ non linear functional defined on $U.$
Let $x\in S$ and $f'(x)\neq 0.$ If $x$ is a critical point of
$f|_S,$ then there exists non zero vector $e_*(x)\in E$ such that\\

 \quad \quad \quad \quad $N(g'(x))\bot e_*(x)$ \quad and \quad $N(f'(x))\bot e_*(x).$
\end{th}
{\it Proof.}\quad Since $f'(x)\neq 0,$ it is clear that
$Nf'(x))\nsubseteq E.$ Hence there exists a non zero vector
$e_*(x)\in E$ such that $e_*(x)\bot Nf'(x)).$ By Theorem 2.1 ,
$e_*(x)\bot N(g'(x)).$
\hfill $\Box$\\

 The following examples
 illustrate an application and significance of Theorem 2.1, although all of them are
 very simple. \\

\begin{ex}\label{1} Suppose that $S$ is the  unit circle in
${\mathbb{R}}^2, x^2+y^2=1,  r_0=\sqrt{{x_0}^2+{y_0}^2}>0,$ and
$f(x,y)=(x-x_0)^2+(y-y_0)^2.$  By Theorem 2.1 go to find the extreme
point $(x,y)$ of $f|_S.$
\end{ex}
It is clear that 1 is a regular value of $g(x,y)=x^2+y^2.$ By
computing simple $$N(f'(x,y))=\{(\Delta x, \Delta y): (x-x_0)\Delta
x +(y-y_0)\Delta y=0\}$$ and $$N(g'(x,y))=\{(\Delta x, \Delta y):
x\Delta x+ y\Delta y =0 \}.$$ If $(x_0, y_0)$ is in $S,$ then
directly, $$f(x_0, y_0)=0\leq (x-x_0)^2+(y-y_0)^2=f(x,y)$$ for any
$(x,y)\in S,$ and so, $(x_0, y_0)$ is an extreme point of $f|_S.$
Otherwise, by Theorem2.1 and
$$dim N(f'(x, y))= dim N(g'(x, y))=1 \quad \quad \forall (x, y)\in
S.$$  it follows that there exists a real number $t$ such that
$$x-x_0=tx \quad y-y_0=ty.$$ (Here $t$ is just Lagrange multiplier.)
So,\\

\quad \quad \quad \quad \quad \quad  $(1-t)x=x_0,$ \quad and \quad
$(1-t)=y_0.$\\

In addition, $(x, y)\in S.$ We then conclude $(1-t)^2={r_0}^2.$
Finally, we get\\

\quad \quad \quad \quad \quad \quad $(x,y)=\frac{1}{r_0}(x_0,y_0)$
\quad or \quad  $-\frac{1}{r_0}(x_0,y_0).$\\

 It is easy to examine
that $(x,y)=\frac{1}{r_0}(x_0,y_0)$ is the required
extreme point.\\

The next example is for generalized regular constraint but not for
regular.

\begin{ex}\label{2} Define $g:{\mathbb{R}}^3\rightarrow
{\mathbb{R}}^3$ as follows
$$g(x)=({x_1}^2+{x_2}^2+{x_3}^2, x_3, x_3) \quad \forall x=(x_1,x_2,x_3)\in {\mathbb{R}}^3. $$
Let $y_0=(1,0,0),$  $S=g^{-1}(y_0)=\{(x_1,x_2,0):
{x_1}^2+{x_2}^2=1\},$ and
$f(x)=(x_1-x_1^0)^2+(x_2-x_2^0)^2+(x_3-x_3^0)^2$ where
$(x_1^0)^2+(x_2^0)^2>0.$ By Theorem 2.1 go to find the extreme point
of $f|_S.$
\end {ex}
It is easy to observe   $$N(f'(x))=\{(\Delta x_1,\Delta x_2, \Delta
x_3): (x_1-x_1^0)\Delta x_1+(x_2-x_2^0)\Delta x_2 -x_3^0\Delta
x_3=0\}$$ and $$N(g'(x))=\{(\Delta x_1, \Delta x_2,0): x_1\Delta
x_1+x_2 \Delta x_2=0\},$$ for $x\in S.$  Obviously,
 $R(g'(x))=\{(2x_1\Delta x_1+2x_2\Delta x_2, \Delta x_3, \Delta x_3 ):
 \forall (\Delta x_1, \Delta x_2, \Delta x_3 )\in {\mathbb{R}}^3\}$ for any $x\in S,$
  so that $Rank(g'(x))=2 $ for all $x\in S,$ Hence, $S$ is the generalized
  regular constraint as indicated in Section 1, but not regular because of $Rank(g'(x))<3.$
 Now we are going to find the extreme point of $f|_S$ by using Theorem
 2.1, it follows from $N(f'(x))\supset N(g'(x))$
 $$N_0=\{(\Delta x_1, \Delta x_2,0) : (x_1-{x_1}^0)\Delta x_1+(x_2-{x_2}^0)\Delta x_2=0\}
 \supset N(g'(x))$$ for any $ x\in S.$
  If $(x_1^0, x_2^0,0)\in S$ then
 directly, $$f(x_1^0, x_2^0,0)=(x_3^0)^2 \leq
 (x_1-x_1^0)^2+(x_2-x_2^0)^3+(x_3^0)^2=f(x)$$ for any $x\in S,$ and
 so, $(x_1^0, x_2^0,0)$ is an extreme point of $f|_S.$ Otherwise, $dim
 N_0=dim N(g'(x))=1,$ so that $N_0=N(g'(x)).$ Thus, similar to
 Example 1, there exists a real number $t$ such that\\

 \quad \quad \quad \quad \quad \quad $x_1-x_1^0=tx_1$ \quad and \quad $x_2-x_2^0=tx_2.$\\

  Let $r_0=\sqrt{(x_1^0)^2+(x_2^0)^2},$
  we can get\\

  \quad \quad \quad \quad \quad $x_1=\frac{1}{r_0} x_1^0$ \quad and \quad $x_2=\frac {1}{r_0}
  x_2^0.$\\

   We now conclude that $(\frac{1}{r_0 }x_1^0,\frac{1}{r_0}x_2^0, 0)$ is an extreme point of
  $f|_S.$ \\
\begin{ex}\label{3} Let $g$ and $f$ be as the same as Example 2.
Apply Theorem 2.2 to find the extreme point of $f|_S.$
\end{ex}
As indicated in Example 2, $$N(f'(x))=\{(\Delta x_1,\Delta x_2,
\Delta x_3): (x_1-x_1^0)\Delta x_1+(x_2-x_2^0)\Delta x_2
-x_3^0\Delta x_3=0\}$$ and $$N(g'(x))=\{(\Delta x_1, \Delta x_2,0):
x_1\Delta x_1+x_2 \Delta x_2=0\},$$ for $x\in S.$ Obviously,
$e_*(x)= (x_1-x_1^0)e_1+(x_2-x_2^0)e_2-x_3^0e_3,$ where
$\{e_i\}_{i=1}^3$ is an orthonormal basis in ${\mathbb{\mathbb{R}}
}^3.$ By Theorem 2.2 one concludes that if $x$ is a critical point
of $f'|_S,$ then
$$(x_1-x_1^0)\Delta x_1 +(x_2-x_2^0)\Delta x_2= -x_1^0\Delta x_1
-x_2^0\Delta x_2=0$$ for all $(\Delta x_1, \Delta
 x_2, 0) \in N(g'(x)).$ Hence $x_1=\lambda x_1^0 $ and $x_2=\lambda
 x_2^0$ for some real number $\lambda .$ Finally, It follows from
$(x_1, x_2, 0)\in S $ that $x_1= \frac{1}{r_0} x_1^0$ and
$x_2=\frac{1}{r_0} x_2^0 $ where $r_0=\sqrt{(x_1^0)^2+(x_2^0)^2},$

   These results of Examples 1 and 2 are very intuitive. The next
   example is not for finding the  critical point, but for showing
   the
   significance of Theorem 2.1.

\begin{ex}\label{4} The constraint $S$ and the non linear function $f(x,y)$ are defined
by\\

\quad \quad $\frac{x^2}{a^2} + \frac{y^2}{b^2}=1$ (both of $a$ and
$b>
0$) \quad and \quad $(x-x_0)^2+(y-y_0)^2 ,$\\

respectively. Let the parameter equation of $S$ be as follows
$$x=a \ cos\theta \quad \rm{and} \quad  y=b\  sin\theta , \ \ 0\leq \theta
<2\pi ,$$ It is known well that if $(a cos\theta ,b sin\theta )$ is
an extreme point of $f|_S,$ then

 \begin{equation}\label{2.2}0=\frac{df(a \ cos\theta ,b \ sin\theta
)}{d\theta }=(b^2-a ^2)sin\theta \ cos\theta +x_0 \ a \ sin \theta
-y_0 \ b \ cos \theta .
\end{equation}
Let us obverse what yields from Theorem 2.1.
\end{ex}
If $(x_0, y_0)\in S,$ it is clear $$f(x_0, y_0)=0\leq
(x-x_0)^2+(y-y_0)^2=f(x, y) \quad \quad \forall (x, y)\in S,$$ i.e.,
$(x_0,y_0)$ is an extreme point of $f|_S.$ So in what follows, we
assume that $(x_0, y_0)$ is not in $S.$  Let $g(x,y)=\frac{x^2}{a^2}
+ \frac{y^2}{b^2}.$ It is easy to examine that 1 is regular value of
$g.$ By  direct computing,
$$N(g'(x,y))=\{ (\Delta x, \Delta y):\frac {x}{a^2}\Delta x+ \frac {y}{b^2}\Delta y=0\}$$
 and $$N(f'(x,y))=\{(\Delta x, \Delta y): (x-x_0)\Delta x +
 (y-y_0)\Delta y=0\}.$$ Since $(x_0, y_0)$ is not in $S$ we conclude
 $dim N(g'(x,y))=dim N(f'(x,y)=1.$ Thus, $N(g'(x,y))=N(f'(x,y)$ by
 Theorem2.1.  Hereby we infer $$\frac{(x-x_0)y}{b^2}-\frac
 {(y-y_0)x}{a^2}=0.$$ Replace  $x$ and $y$ in the equation
 above by $a \ cos\theta $ and $b \ sin \theta,$ respectively, then we get
 Equation (2.2).\\

 The example shows that Theorem 2.1 implies the classical principle (2.2).
\begin{re}\label{2.1} \ In case $F=\mathbb{R},$ the classical principle
for seeking the extreme point of $f|_S$ is to solve Euler equations
with Lagrange multiplier $L,$ $f'(x)-L g'(x)=0$ and $g(x)=y_0.$ In
this case, $L$ is a real number. When $dimF>1,$ $L$ is a bounded
linear functional on $F.$ In the sequel, it will be proved that in
the case of regular constraint $S,$ $L$ is unique; otherwise $L$ is
ill-posed. Hence under the generalized regular but not regular
constraint, it is very difficult to solve Euler equations with
Lagrange multiplier. So, recently, researching of critical points
under the non regular constraint has been a concerned focus in
optimization theory. Therefore, Theorems 2.1 and 2.2 seem to be a
new and available way.
\end{re}

\section{Lagrange  Multiplier    }\label{s3}
 \vspace{0.3cm}
Let $A\in B(E,F)$ be double split. As indicated in Section 1, there
exists a generalized inverse $A^+ \in B(F,E)$ of $A$ such that
$A^+AA^+=A^+$ and $AA^+A=A.$ It is easy to obverse
$A^+A=P_{R(A^+)}^{N(A)}$ and $AA^+=P_{R(A)}^{N(A^+)},$ where
$P_{R(A^+)}^{N(A)}$ is a projection from $E$ onto $R(A^+)$ and
$P_{R(A)}^{N(A^+)}$ from $F$  onto $R(A)$ coordinate to the
following decompositions : $E=N(A)\oplus R(A^+)$ and $F=R(A)\oplus
N(A^+),$ respectively. So, $I_E-
P_{R(A^+)}^{N(A)}=P_{N(A)}^{R(A^+)}$ and
$I_F-P_{R(A)}^{N(A^+)}=P_{N(A^+)}^{R(A)}.$ In this section, we
discuss the classical critical point principle with Lagrange
multiplier $L$ under generalized regular constraint,  give an
express of $L$ by using generalized inverse and show that $L$ is
unique under regular constraint.
\begin{th}\label{3.1} Suppose that $g :U\subset E\rightarrow F$ is
a $c^1$ map from open set $U$ in Banach space $E$ into Banach space
$F.$ Let $g'^+(x)$ be a generalized inverse of $g'(x), f $ a $ c^1$
non linear functional on $U,$ and $S=$ the preimage $g^{-1}(y_0),
y_0\in F.$ If $y_0$ is a generalized regular value of $g,$ and $x\in
U$ is a critical point of $f|_S,$ then
$$f'(x)-f'(x)\circ {g'}^+(x)\circ g'(x)=0
\quad {\rm{and}} \quad g(x)=y_0,$$ for any ${g'}^+(x) \in
GI(g'(x)).$
 i,e., $f'(x)\circ {g'}^+(x)\in F^*$ is a Lagrange multiplier.
\end{th}
{\it Proof.}\quad  Since
$$R(I_E-{g'}^+(x)g'(x))=R(P_{N(g'(x))}^{R(g'^+(x)})=N(g'(x))$$ and
by Theorem 2.1,$$0=f'(x)(I_E-{g'}^+(x)g'(x))=f'(x)-f'(x)\circ
{g'}^+(x)\circ g'(x)$$ whenever $x$ is a critical point of $f|_S.$
This
proves the theorem.\hfill $\Box$\\
In what follows, we are going to show that Lagrange multiplier is
unique under regular constraint.\\
\begin{th}\label{3.1} If $S=$ the preimage $g^{-1}(y_0)$ is the regular
constraint, i.e., $y_0\in F$ is regular value of $g,$ then the
Lagrange multiplier $L$ is unique.
\end{th}
{\it Proof.}\quad  Assume $x\in S$ is a critical point of $f|_S.$
Let $q_1(x)$ and $q_2$ be arbitrary two generalized inverses of
$g'(x).$
 Since $y_0$ is a regular value, $R(g'(x))=F.$ Hence.
 $$(f'(x)\circ q_i(x))(y)=(f'(x)\circ q_i(x)\circ g'(x))(e)=f'(x)e,  \quad i=1,2, \quad \forall y\in F,$$
 where  $y=g'(x)(e).$ Therefore $f'(x)\circ q_1(x)=f'(x)\circ q_2(x).$ Let
 $L$ be the Lagrange multiplier. Similarly,
$$ L(y)=(L\circ g'(x))(e)=f'(x)(e)=(f'(x)\circ g'^+(x))(y) \quad
 \forall y\in F.$$ The theorem is proved. $\Box$\\
\section{ Differential Construction of $GI(A)$ and Ill-Posed Lagrange Multipliers }\label{s1}
 \vspace{0.3cm}
  In this section, we are going to discuss the differential construction
  of $GI(A)$ and prove that $GI(A)$ is smooth diffeomorphic to some
  Banach space. Then we show that when $S$ is a generalized regular
  but not regular
  constraint,
   the Lagrange multiplier $L$ is
  ill-posed.\\

  Let $A\in B(E,F)$ be double split, and $A_0^+\in GI(A)$ coordinating
  to the decompositions : $E=N(A)\oplus R_0^+$ and $F=R(A)\oplus
  N_0^+,$ i.e., $R_0^+=R(A_0^+)$ and $N_0^+=R(I_F-AA_0^+).$
  Consider the following map $M:B(R_0^+,N(A))\times B(N_0^+,R(A)))\rightarrow
  B(F,E),$
.$$M(\alpha , \beta )=(I_E+\alpha)A_0^+(I_F-\beta P_{N_0^+}^{R(A)})
                  \quad \forall \ (\alpha ,\beta)\in B(R_0^+,N(A))\times
B(N_0^+,R(A))).$$

\begin{lem}\label{4.1} Let $B=M(\alpha , \beta ) \  \forall(\alpha  , \beta )\in
B(R_0^+, N(A)) \times B(N_0^+, R(A)).$  The following preperties of
$M$
for any $(\alpha ,\beta ) \in B(R_0^+, N(A)) \times B(N_0^+, R(A))$ hold:\\
  $$R(B)=\{e+\alpha (e): \forall e\in R_0^+\},$$
$$ N(B)=\{d +\beta (d): \forall d\in N_0^+\},$$ and $M(\alpha ,
\beta )\in GI(A).$

\end{lem}
{\it Proof.}\quad  Evidently,
 $$
BA=(I_E+\alpha )A_0^+(I_F-\beta P_{N_0^+}^{R(A)})A=(I_E+\alpha
)A_0^+A, \ {\rm{and}} \ AB=AA_0^+(I_F-\beta P_{N_0^+}^{R(A)}) \
$$
Hereby, $ABA=A(I_E+\alpha )A_0^+A=AA_0^+A=A ,$and $ BAB=(I_e+\alpha
)A_0^+AA_0^+(I_F-\beta P_{N_0^+}^{R(A)})=B. $ This says $B\in
GI(A).$ Note that $e+\alpha (e)=0 $ for $e\in R_0^+$ implies $e=0.$
Then we have $$y\in N(B) \Leftrightarrow A_0^+(I_F-\beta
P_{N_0^+}^{R(A)})y=0 \Leftrightarrow (I_F-\beta P{N_0^+}^{R(A)})y\in
N_0^+ $$
$$\Leftrightarrow P_{N_0^+}^{R(A)}(I_F-\beta P_{N_0^+}^{R(A)})y=(I_F-\beta P_{N_0^+}^{R(A)})y\Leftrightarrow
 P_{N_0^+}^{R(A)}y=(I_F-\beta P_{N_0^+}^{R(A)})y.$$ Hence  $y=
 P_{N_0^+}^{R(A)}y +\beta (P_{N_0^+}^{R(A)}y) \ i,e,. N(B)\subset \{d+\beta (d) : \forall d\in N_0^+\} .$
 Conversely, let $y=d+\beta (d)$ for any $d\in N_0^+,$ then
 $P_{N_0^+}^{R(A)}y=d $ and so $y=P_{N_0^+}^{R(A)}y +\beta
 (P_{N_0^+}^{R(A)}y).$ Thus, by the equivalent relations above, $y\in
 N(B).$ This shows $N(B)=\{d +\beta (d): \forall d\in N_0^+\}.$ Let
 $x\in R(B),$ $x=(I_E+\alpha )A_0^+(I_F-\beta
  P_{N_0^+}^{R(A)})y,$ and $e=A_0^+(I_F-\beta
   P_{N_0^+}^{R(A)})y.$ Then $x=e+\alpha (e)$ and so, we prove $R(B)\subset
 \{e+\alpha (e): \forall e\in R_0^+\}.$ Conversely, let $x=e+\alpha
 (e)$ where $e\in R_0^+,$ and $e=A_0^+y, \ y\in R(A).$ Then
 $$x=(I_E+\alpha )A_0^+y=(I_E+\alpha )A_0^+ (I_F-\beta P_{N_0^+}^{R(A)})y$$
since $P_{N_0^+}^{R(A)}y=0.$ Now, one concludes $R(B)=\{e+\alpha
(e): \forall e\in R_0^+\}.$
$\Box$\\

\begin{lem}\label{4.2} $ M(\alpha , \beta ): B(R_0^+,N(A))\times B(N_0^+,R(A))\rightarrow
 GI(A)$ is bijective.
\end{lem}

{\it Proof.}\quad  It is easy to see that if $B,B_1\in GI(A),$
$R(B)=R(B_1)$ and  $N(B)=N(B_1),$ then $B=B_1.$ In fact, $$B=
BAB=BP_{R(A)}^{N(B)}=BP_{R(A)}^{N(B_1)} =BAB_1
$$ $$=P_{R(B)}^{N(A)}B_1=P_{R(B_1)}^{N(A)}B_1=B_1AB_1=B_1.$$
Hereby, due to the conclusions about $R(B)$ and $N(B)$ in Lemma 4.1
 one can infer that $M(\alpha ,\beta )$ is injective. Let both of $B$ and
 $A_0^+$ be in $GI(A).$ we have
$$ P_{R_0^+}^{N(A)}P_{R(B)}^{N(A)}e=
P_{R_0^+}^{N(A)}(P_{R(B)}^{N(A)}e+P_{N(A)}^{R(B)}e)=P_{R_0^+}^{N(A)}e=e
\ \forall e\in R_0^+$$ and
$$P_{R(B)}^{N(A)}P_{R_0^+}^{N(A)}e=P_{R(B)}^{N(A)}(P_{R_0^+}^{N(A)}e+P_{N(A)}^{R_0^+}e)=P_{R(B)}^{N(A)}e
=e \ \forall e\in R(B).$$ Similarly,\\

\quad \quad \quad \quad $P_{N_0^+}^{R(A)}P_{N(B)}^{R(A)}d=d \
\forall d\in N_0^+$ \quad and\quad
$P_{N(B)}^{R(A)}P_{N_0^+}^{R(A)}d=d
 \ \forall d\in N(B).$\\

  Let\\

 \quad \quad \quad \quad \quad \quad $\alpha
 =P_{N(A)}^{R_0^+}P_{R(B)}^{N(A)}$
 \quad and \quad $\beta =P_{R(A)}^{N_0^+} $\\

  for any $B\in
 GI(A).$ It is easy to check\\

 \quad \quad $R(B)=\{e+\alpha (e): \forall e\in R_0^+\}
\quad $ and $\quad N(B)=\{d+\beta (d): \forall e\in N_0^+\}.$\\

  In fact,
 $$ h=P_{R_0^+}^{N(A)}h+ P^{R_0^+}_{N(A)}h=P_{R_0^+}^{N(A)}h+ P^{R_0^+}_{N(A)}
 P_{R(B)}^{N(A)}P_{R_0^+}^{N(A)}h=P_{R_0^+}^{N(A)}h+\alpha(P_{R_0^+}^{N(A)}h)
 $$ for any $h\in R(B),$ so that $R(B)\subset \{e+\alpha (e): \forall
 e\in R_0^+\};$ conversely,
 $$h=e+\alpha
 (e)=P_{R_0^+}^{N(A)}P_{R(B)}^{N(A)}e+P_{N(A)}^{R_0^+}P_{R(B)}^{N(A)}e=P_{R(B)}^{N(A)}e\in
 R(B)$$
  for any$ e\in R_0^+.$ Similarly, one can verify the above relation
  about $N(B),$ By Lemma 4.1 it follows $B=M(\alpha , \beta )$
  since both of $B$ and $M(\alpha ,\beta )$ are in $GI(A).$
  Therefore,$ M:B(R_0^+, N(A))\times B(N_0^+, R(A))\rightarrow GI(A)$
  is bijective.
$\Box$\\
\begin{th}\label{4.1} The mapping $M(\alpha , \beta )$ from  $B(R_0^+,N(A))\times
B(N_0^+, R(A))$ onto $GI(A)$ is smooth diffeomorphism.
\end{th}

{\it Proof.}\quad By direct computing  \\

$M(\alpha +\Delta \alpha ,\beta +\Delta \beta )$\\

$=(I_E+\alpha )A_0^+(I_F-(\beta +\Delta \beta
)P_{N_0^+}^{R(A)}+\Delta \alpha A_0^+(I_F-(\beta +\Delta \beta
)P_{N_0^+}^{R(A)})$\\

 $=M(\alpha ,\beta )-(I_E+\alpha )A_0^+\Delta
\beta P_{N_0^+}^{R(A)}+\Delta \alpha A_0^+(I_F-\beta
P_{N_0^+}^{R(A)})-\Delta \alpha A_0^+\Delta \beta
P_{N_0^+}^{R(A)}.$\\

 Hereby we infer\\

  $DM(\alpha , \beta )< \Delta
\alpha , \Delta \beta
>$\\

$=\Delta \alpha A_0^+(I_F-\beta P_{N_0^+}^{R(A)})-(I_E+\alpha
)A_0^+\Delta
\beta P_{N_0^+}^{R(A)}.$\\

Similarly,\\

$DM(\alpha +\Delta {\alpha } _1, \beta +\Delta {\beta}_1)<\Delta
\alpha
,\Delta \beta >$,\\

$=\Delta \alpha A_0^+(I_F-(\beta +\Delta
{\beta}_1)P_{N_0^+}^{R(A)}-(I_E+(\alpha+\Delta {\alpha
}_1))A_0^+\Delta \beta P_{N_0^+}^{R(A)}$\\

$=\Delta \alpha A_0^+(I_F-\beta )P_{N_0^+}^{R(A)}-\Delta \alpha
A_0^+\Delta {\beta }_1P_{N_0^+}^{R(A)} - (I_E+\alpha )A_0^+\Delta
\beta P_{N_0^+}^{R(A)}-\Delta {\alpha }_1A_0^+\Delta \beta
P_{N_0^+}^{R(A)}$\\

$=DM(\alpha ,\beta )<\Delta \alpha , \Delta \beta>-\Delta {\alpha
}A_0^+\Delta {\beta}_1P_{N_0^+}^{R(A)}-\Delta {\alpha}_1 A_0^+\Delta
\beta P_{N_0^+}^{R(A)}.$\\

Hereby we conclude\\

 $D^2M(\alpha , \beta )< (\Delta \alpha , \Delta
\beta ), (\Delta {\alpha }_1, \Delta {\beta}_1)
>$ \\

$=-\Delta \alpha A_0^+\Delta {\beta }_1 P_{N_0^+}^{R(A)}- \Delta
{\alpha }_1
 A_0^+\Delta \beta P_{N_0^+}^{R(A)}.$\\

 Note that $D^2M (\alpha , \beta )$ is independent of $\alpha $ and
 $\beta .$ Hence $$D^n(\alpha ,\beta )=0 \quad n\geq 3.$$ The
 theorem is proved due to Lemma 4.2.
$\Box$\\
For the details see \cite {m4} and \cite{m5}.\\
 We are now in the
position to consider the ill-posed Lagrange multiplier in Euler
equations under  non regular constraint.

\begin{th}\label{4.2} Suppose that $g:U\subset E\rightarrow F$ is
a $c^1$ map, and $S=$ the preimage $g^{-1}(y_0)$ is a generalized
regular constraint but not regular, i,e,. $y_0\in F$ is a
generalized regular value of $g$ but not regular. Let $f$ be a non
linear functional defined on $U.$ Then the Lagrange multiplier $L$
in Euler equations is ill-posed.
\end {th}

{\it Proof.}\quad Let $x$ be a critical point of $f|_S,$ and $e_0\in
E$  such that $f'(x)e_0=1.$ Since $x$ is a critical point of $f|_S$
by Theorem 3.1,
$$1=(f'(x)\circ {g'}^+(x)\circ g'(x))(e_0) \quad  {\rm{for}} \  {\rm{any}} \
{\rm{fixed}} \ {g'}^+(x)\in GI(g'(x)).$$ Let $y_1 =g'(x)e_0 \in
R(g'(x)).$ Since $x\in S$ and $S$ is the non regular constraint,
there exists a non zero $y^+$ in $N({g'}^+(x)).$  Write
$N({g'}^+(x))$ by $N_0^+,$
 and
 let $N_0^+=[y^+] \oplus N^+,$ where $[y^+]$ denotes
the subspace generated by $y^+.$ Define $\beta \in B(N_0^+,$
$R(g'(x))$ as follows $$ \beta (y)=0 \quad {\rm{if}} \quad  y\in
N^+; \quad {\rm{else}} \quad \beta (y^+)=y_1.$$ By Lemma 4.1, $$
B=g'(x)(I_F-\beta P_{N_0^+}^{R(g'(x))}) \in GI(g'(x)) \quad
{\rm{and}} \quad y^+ +y_1 \in N(B). $$ Now we see that two deferent
Lagrange multipliers at the critical point $x,$ $L=f'(x)\circ
{g'}^+(x)$ and $L_1=f'(x)\circ B,$ both  are in $F^*$ and satisfy
$$ L=(f'(x)\circ {g'}^+(x))(y^+ +y_1)=(f'(x)\circ
{g'}^+(x))(y_1)=((f'(x)\circ {g'}^+(x))(g'(x))e_0)=1,$$ and
$$L_1=(f'(x)\circ B)((y^+ +y_1)=0.$$
The theorem is proved.
$\Box$\\

The constraint $S$ in Example 2 is non regular. Now take it with
$x_3^0\neq 0$ for example to illustrate the ill-posed Lagrange
multiplier in case of non regular constraint.

Let $ \{e_i\}_1^3$ and $ \{\varepsilon _i\}_1^3$ denote orthonormal
bases of  $ {\mathbb{R} }^3 $  containing the domain and the range
of $g,$ respectively. As is indicated in Example 2 that
$x_0=\frac{1}{r_0}(x_1^0, x_2^0, 0)$ is a critical point of $f|_S.$
By direct computing,
$$g'(x_0)(x_1e_1+x_2e_2+x_3e_3)=\frac{1}{r_0}(x_1^0x_1+x_2^0x_2)\varepsilon
_1+x_3\varepsilon _2+x_3\varepsilon _3$$ and
$$f'(x_0)(y_1\varepsilon _1+y_2\varepsilon _2+y_3\varepsilon
_3)=2(\frac{2}{r_0} -1)y_1+2(\frac{2}{r_0}
-1)y_2-2x_3^0y_3.$$Hereby,
$$N(g'(x_0))=\{x_1e_1+x_2e_2: x_1^0x_1+x_2^0x_2=0\},$$ and so,
$$N(g'(x_0))\bot \ (x_1^0e_1+x_2^0e_2) \  \bot x_3^0e_3.$$ Then, we have the following
decomposition of  ${\mathbb{R} }^3$ containing the domain of $g$
$${\mathbb{R} }^3 =N(g'(x_0)\oplus [x_1^0e_1+x_2^0e_2] \oplus
[e_3],$$ where $[,]$ denotes the subspace generated by the vector in
the bracket, and $\oplus $ orthogonal sum.
 Moreover, since
$g'(x_0)(x_1^0e_1+x_2^0e_2)= 2r_0\varepsilon _1$ and
$g'(x_0)(x_3^0e_3)=x_3^0\varepsilon _2+x_3^0\varepsilon _3$ it
follows
$$R(g'(x_0))=[2r_0\varepsilon _1 ]\oplus [x_3^0\varepsilon
_2+x_3^0\varepsilon _3].$$ Hereby we infer the decomposition of
${\mathbb{R} }^3$ containing the range of $g$:$${\mathbb{R}
}^3=R(g'(x_0))\oplus [-x_3^0\varepsilon _2+x_3^0\varepsilon _3].
$$ Let $N_0^+=[-x_3^0\varepsilon _2+x_3^0\varepsilon _3] .$ Now we can define ${g'(x_0)}^+ \in GI(g'(x_0))$ as follows,
$$ {g'(x_0)}^+ y=0 \quad   {\rm{if}} \quad   y\in N_0^+,$$
$$ \quad \quad \quad \quad \quad \quad \quad \quad   =x_1^0e_1+x_2^0e_2 \quad   {\rm{if}} \quad y=2r_0\varepsilon _1,$$
$$ \quad \quad \quad \quad \quad \quad \quad \quad   =x_3^0e_3 \quad {\rm{if}} \quad y=x_3^0\varepsilon _2 +x_3^0\varepsilon _3 .$$
Obviously, $L=f'(x_0)\circ {g'(x_0)}^+ \in F^*$ and\\

$ L(2x_3^0\varepsilon _3)=(f'(x_0)\circ
{g'(x_0)}^+)((-x_3^0\varepsilon
_2+x_3^0\varepsilon _3)+(x_3^0\varepsilon _2+x_3^0\varepsilon _3))$\\

$ \quad \quad \quad \quad =(f'(x_0)\circ
{g'(x_0)}^+)(x_3^0\varepsilon
_2+x_3^0\varepsilon _3)$\\

$ \quad \quad \quad \quad =f'(x_0)(x_3^0 e_3)=-2(x_3^0)^2\neq 0.$\\

Let $\beta \in B(N_0^+, R(g'(x_0)))$ such that
$$ \beta (-x_3^0\varepsilon
_2+x_3^0\varepsilon _3)= x_3^0\varepsilon _2+x_3^0\varepsilon _3.$$

Set $$B={g'(x_0)}^+(I_F-\beta P_{N_0^+}^{R(g'(x_0))}).$$

By Lemma 4.1, $$2x_3^0\varepsilon _3=(-x_3^0\varepsilon
_2+x_3^0\varepsilon _3)+(x_3^0\varepsilon _2+x_3^0\varepsilon _3)\in
N(B).$$ (This conclusion can also be to check directly.) Hence
$$L_1 (2x_3^0\varepsilon _3)=(f'(x_0)\circ B)(2x_3^0\varepsilon _3)=0.$$
Therefore $L$ and $L_1$ are two deferent Lagrange multipliers.

\begin{re}\label{4.1} \ It is not enough for the generalized preimage theorem
to express several constraints in  optimization theory, finance
mathematics and so on. A generalization to Thom's famous result for
transversality,  generalized transversality theorem is helpful, (see
\cite {m2}). We can also have similar principle to Theorem 2.1,
which will be discussed else where.

\end{re}

1, Tseng Yoan Rong Functional Research Center. Harbin Normal
University, Harbin 150080, P.R.China

2 Department of Mathematics ,Nanjing University, Nanjing, 210093,
P.R.China\\
 E-mail address: \ jipuma@126.com


\begin{thebibliography}{99}

\bibitem [AMR] {abr}\   Abraham R,  Marsden J. E, Ratiu T. \  {\it
Manifold, Tensor analysis and its Applications,} WPC:
Springer-Verlag, (1988).


\bibitem [B] {b}  Berger M, {\it  Nonlinearity and Functional Analysis},
New York Academic Press, (1976).

\bibitem [HM] {hm}\ Qianlian  H, Jipu Ma, {\it Perturbation Analysis
of Generalized Inverses of Linear Operators in Banach Space}, Linear
Algebra Appl. 389 (2004),335-364.


\bibitem [M1]  {m1}\  Jipu  Ma , {\it Three classes of smooth Banach
manifolds in B(E,F)}, Sci.China Ser.A 50:9(2007), 1233-1239.


\bibitem [M2] {m2}\ Jipu Ma , {\it A generalized transversility in
 global analysis}, PJM 236:2(2008), 357-371.


\bibitem [M3] {m3} \ Jipu Ma , { \it Complete rank theorem of
advanced calculus and singularities of bounded linear operators},
Front.Math. china 3:2 (2008), 304-316.


\bibitem [M4] {m4} \ Zhaofeng Ma ,Jipu  Ma , {\it The Smooth Banach
Submanifold $B^*(E,F)$ in $B(E,F)$},  Sci.China Ser.A 52:11(2009),
2479-2492.


\bibitem [M5] {m5} \ Zhaofeng Ma , Jipu Ma , {\i A Common Property of
$R(E,F)$ and ${\mathbb{R}}^n,{\mathbb{R}}^m)$ and A New Method for
Seeking A Path to Connect Two Operators},  Sci.China Ser.A
53:10(2010), 2605-2620.


\bibitem [M6
]  {m6}\  Jipu Ma , {\it A generalized preimage theorem in global
nalysis},  Science in China, Ser.A, 44:(2001), 299-303.


\bibitem [M7] {m7} \ Jipu Ma , {\it (1,2)-Inverse of Operators
between Banach Spaces and Local Conjugacy Theorem}, Chi. Ann. Math..
 20(B):1 (1999), 57-62.


\bibitem [M8] {m8} \ Jipu Ma , {\it Local Conjugacy Theorem, Rank
Theorems in Advanced Calculus and  A Generalized Principle for
Constructing Banagh Manifold }, Science in China, Ser.A, 43:12
(2000), 1233-1237.


\bibitem [M9] {m9} \ Jipu Ma, {\it Rank Theorem of Operators Between
Banach Spaces }, Science in China, Ser.A, 43:1(2000), 1-5.

\bibitem [Z]  {zei}\  Zeidler E,  {\it Nonlinear Functional Analysis and
Its Applications IV}, New York-Berlin: Springer-Verlag, (1988).


\bibitem[N] {n}\  Nashed M. Z, {\it Generalized inverses and
applications,} Academic Press, New York, San-Francisco, London 1976.
\end{thebibliography}
\end{document}